\documentclass{gtart}     

\def\ifplaintex{\expandafter\ifx\csname documentclass\endcsname\relax}

\def\gtp{{\mathsurround=0pt\it $\cal G\mskip-2mu$eometry \&\ 
$\cal T\!\!$opology $\cal P\!$ublications}}  

\def\recd{{\small Received:\qua\receiveddate\ifx\reviseddate\relax
\else\qquad Revised:\qua\reviseddate\fi\par}} 

\def\lognumber#1{\def\thelognumber{#1}}
\def\volumenumber#1{\def\thevolumenumber{#1}}
\def\volumeyear#1{\def\thevolumeyear{#1}}
\def\papernumber#1{\def\thepapernumber{#1}}
\def\pagenumbers#1#2{\def\startpage{#1}\def\finishpage{#2}}
\def\published#1{\def\publishdate{#1}}

\def\received#1{\def\receiveddate{#1}}

\def\accepted#1{\def\accepteddate{#1}}

\long\def\asciiabstract#1{\long\def\theasciiabstract{#1}}

\let\\\par\let\thelognumber\relax\let\thevolumenumber\relax
\let\thepapernumber\relax\let\thevolumeyear\relax\let\startpage\relax
\let\finishpage\relax\let\publishdate\relax\let\receiveddate\relax
\let\reviseddate\relax\let\accepteddate\relax\let\theasciititle\relax
\let\theasciiauthors\relax
\let\theasciiabstract\relax

\let\theasciiemail\relax

\ifplaintex
\font\logobig=cmssbx10 scaled 3836
\font\logomed=cmssbx10 scaled 2557
\else
\font\logobig=cmssbx10 scaled 4200
\font\logomed=cmssbx10 scaled 2800
\fi

\long\def\makeagttitle{   
\count0=\startpage
\agt\hfill      
\hbox to 45truept{\vbox to 0pt{\vglue -13truept{\logomed A\kern -.37em{\logobig 
T}\kern -.38em G}\vss}\hss}
\break
{\small Volume \thevolumenumber\ (\thevolumeyear)
\startpage--\finishpage\nl
Published: \publishdate}

\vglue .25truein

{\parskip=0pt\leftskip 0pt plus
1fil\def\\{\par\smallskip}{\Large\bf\thetitle}\par\medskip} \vglue
0.05truein

{\parskip=0pt\leftskip 0pt plus 1fil\def\\{\par}{\sc\theauthors}
\par\medskip}%
 
\vglue 0.03truein 

{\small\leftskip 25truept\rightskip 25truept{\bf Abstract}\stdspace\theabstract

{\bf AMS Classification}\stdspace\theprimaryclass
\ifx\thesecondaryclass\relax\else; \thesecondaryclass\fi\par
{\bf Keywords}\stdspace \thekeywords\par}\vglue 7truept

}   

\ifplaintex
\font\phead=cmsl9 scaled 950
\font\pnum=cmbx10 scaled 913
\font\pfoot=cmsl9 scaled 950
\headline{\vbox to 0pt{\vskip -4.5mm\line{\small\phead\ifnum
\count0=\startpage ISSN 1472-2739 (on-line) 1472-2747 (printed)
\hfill {\pnum\folio}\else\ifodd\count0\def\\{ }%
\ifx\theshorttitle\relax\thetitle\else\theshorttitle\fi\hfill{\pnum\folio}
\else\def\\{ and }{\pnum\folio}\hfill\ifx\theshortauthors\relax\theauthors
\else\theshortauthors\fi\fi\fi}\vss}}
\footline{\vbox to 0pt{\vglue 0mm\line{\small\pfoot\ifnum\count0=\startpage
\copyright\ \gtp\hfill\else
\agt, Volume \thevolumenumber\ (\thevolumeyear)\hfill\fi}\vss}}
\else
\font=cmsl9 scaled 1050
\font\lnum=cmbx10 
\font=cmsl9 scaled 1050
\makeatletter
\def\@oddhead{{\small\lhead\ifnum\count0=\startpage ISSN 1472-2739 
(on-line) 1472-2747 (printed)\hfill {\lnum\number\count0}\else\ifodd\count0
\def\\{ }\ifx\theshorttitle\relax \thetitle \else\theshorttitle\fi\hfill
{\lnum\number\count0}\else\def\\{ and }{\lnum\number\count0}
\hfill\ifx\theshortauthors\relax 
\theauthors\else\theshortauthors\fi\fi\fi}}\def\@evenhead{\@oddhead}
\def\@oddfoot{\small\lfoot\ifnum\count0=\startpage\copyright\ \gtp\hfill\else
\agt, Volume \thevolumenumber\ (\thevolumeyear)\hfill\fi}
\def\@evenfoot{\@oddfoot}
\makeatother
\fi
\let\maketitlepage\makeagttitle
\let\makeshorttitle\maketitlepage
\let\maketitle\maketitlepage

\newwrite\gtoutfile
\long\gdef\makeheadfile{  
{\def\\{, }\def\s{ }
\immediate\openout\gtoutfile head.xxx
\immediate\write\gtoutfile{To: math@arxiv.org}
\immediate\write\gtoutfile{Subject: put OR rep NNNNN:ppppp}
\immediate\write\gtoutfile{--text follows this line--}
\immediate\write\gtoutfile{Proxy-for: \ifx\theasciiauthors\relax
\theauthors\else\theasciiauthors\fi\s<\ifx\theasciiemail\relax\theemail\else\theasciiemail\fi>}
\immediate\write\gtoutfile{\noexpand\\}
\immediate\write\gtoutfile{Authors: \ifx\theasciiauthors\relax
\theauthors\else\theasciiauthors\fi}
{\def\\{ }\immediate\write\gtoutfile{Title: \ifx\theasciititle\relax
\thetitle\else\theasciititle\fi}}
\immediate\write\gtoutfile{Subj-class: GT or SG, GR etc}
\immediate\write\gtoutfile{MSC-class: \theprimaryclass\ifx\thesecondaryclass\relax\else, \thesecondaryclass\fi}
\immediate\write\gtoutfile{Journal-ref: Algebraic and Geometric Topology \thevolumenumber\s
(\thevolumeyear) \startpage-\finishpage}
\immediate\write\gtoutfile{Comments: Published by Algebraic and
Geometric Topology at}
\immediate\write\gtoutfile{\s\s\s  http://www.maths.warwick.ac.uk/agt/AGTVol\thevolumenumber/agt-\thevolumenumber-\thepapernumber.abs.html}
\immediate\write\gtoutfile{\noexpand\\}
\immediate\write\gtoutfile{}
\ifx\theasciiabstract\relax
\immediate\write\gtoutfile{\theabstract}\else
\immediate\write\gtoutfile{\theasciiabstract}\fi
\immediate\write\gtoutfile{}
\immediate\write\gtoutfile{\noexpand\\}
\immediate\write\gtoutfile{}
\immediate\closeout\gtoutfile}}  

\def\maketitlepage{\makeagttitle\makeheadfile}
\let\makeshorttitle\maketitlepage
\let\maketitle\maketitlepage

\def\ifplaintex{\expandafter\ifx\csname documentclass\endcsname\relax}

\def\gtp{{\mathsurround=0pt\it $\cal G\mskip-2mu$eometry \&\ 
$\cal T\!\!$opology $\cal P\!$ublications}}  

\def\recd{{\small Received:\qua\receiveddate\ifx\reviseddate\relax
\else\qquad Revised:\qua\reviseddate\fi\par}} 

\def\lognumber#1{\def\thelognumber{#1}}
\def\volumenumber#1{\def\thevolumenumber{#1}}
\def\volumeyear#1{\def\thevolumeyear{#1}}
\def\papernumber#1{\def\thepapernumber{#1}}
\def\pagenumbers#1#2{\def\startpage{#1}\def\finishpage{#2}}
\def\published#1{\def\publishdate{#1}}

\def\received#1{\def\receiveddate{#1}}

\def\accepted#1{\def\accepteddate{#1}}

\long\def\asciiabstract#1{\long\def\theasciiabstract{#1}}

\let\\\par\let\thelognumber\relax\let\thevolumenumber\relax
\let\thepapernumber\relax\let\thevolumeyear\relax\let\startpage\relax
\let\finishpage\relax\let\publishdate\relax\let\receiveddate\relax
\let\reviseddate\relax\let\accepteddate\relax\let\theasciititle\relax
\let\theasciiauthors\relax
\let\theasciiabstract\relax

\let\theasciiemail\relax

\ifplaintex
\font\logobig=cmssbx10 scaled 3836
\font\logomed=cmssbx10 scaled 2557
\else
\font\logobig=cmssbx10 scaled 4200
\font\logomed=cmssbx10 scaled 2800
\fi

\long\def\makeagttitle{   
\count0=\startpage
\agt\hfill      
\hbox to 45truept{\vbox to 0pt{\vglue -13truept{\logomed A\kern -.37em{\logobig 
T}\kern -.38em G}\vss}\hss}
\break
{\small Volume \thevolumenumber\ (\thevolumeyear)
\startpage--\finishpage\nl
Published: \publishdate}

\vglue .25truein

{\parskip=0pt\leftskip 0pt plus
1fil\def\\{\par\smallskip}{\Large\bf\thetitle}\par\medskip} \vglue
0.05truein

{\parskip=0pt\leftskip 0pt plus 1fil\def\\{\par}{\sc\theauthors}
\par\medskip}%
 
\vglue 0.03truein 

{\small\leftskip 25truept\rightskip 25truept{\bf Abstract}\stdspace\theabstract

{\bf AMS Classification}\stdspace\theprimaryclass
\ifx\thesecondaryclass\relax\else; \thesecondaryclass\fi\par
{\bf Keywords}\stdspace \thekeywords\par}\vglue 7truept

}   

\ifplaintex
\font\phead=cmsl9 scaled 950
\font\pnum=cmbx10 scaled 913
\font\pfoot=cmsl9 scaled 950
\headline{\vbox to 0pt{\vskip -4.5mm\line{\small\phead\ifnum
\count0=\startpage ISSN 1472-2739 (on-line) 1472-2747 (printed)
\hfill {\pnum\folio}\else\ifodd\count0\def\\{ }%
\ifx\theshorttitle\relax\thetitle\else\theshorttitle\fi\hfill{\pnum\folio}
\else\def\\{ and }{\pnum\folio}\hfill\ifx\theshortauthors\relax\theauthors
\else\theshortauthors\fi\fi\fi}\vss}}
\footline{\vbox to 0pt{\vglue 0mm\line{\small\pfoot\ifnum\count0=\startpage
\copyright\ \gtp\hfill\else
\agt, Volume \thevolumenumber\ (\thevolumeyear)\hfill\fi}\vss}}
\else
\font=cmsl9 scaled 1050
\font\lnum=cmbx10 
\font=cmsl9 scaled 1050
\makeatletter
\def\@oddhead{{\small\lhead\ifnum\count0=\startpage ISSN 1472-2739 
(on-line) 1472-2747 (printed)\hfill {\lnum\number\count0}\else\ifodd\count0
\def\\{ }\ifx\theshorttitle\relax \thetitle \else\theshorttitle\fi\hfill
{\lnum\number\count0}\else\def\\{ and }{\lnum\number\count0}
\hfill\ifx\theshortauthors\relax 
\theauthors\else\theshortauthors\fi\fi\fi}}\def\@evenhead{\@oddhead}
\def\@oddfoot{\small\lfoot\ifnum\count0=\startpage\copyright\ \gtp\hfill\else
\agt, Volume \thevolumenumber\ (\thevolumeyear)\hfill\fi}
\def\@evenfoot{\@oddfoot}
\makeatother
\fi
\let\maketitlepage\makeagttitle
\let\makeshorttitle\maketitlepage
\let\maketitle\maketitlepage

\newwrite\gtoutfile
\long\gdef\makeheadfile{  
{\def\\{, }\def\s{ }
\immediate\openout\gtoutfile head.xxx
\immediate\write\gtoutfile{To: math@arxiv.org}
\immediate\write\gtoutfile{Subject: put OR rep NNNNN:ppppp}
\immediate\write\gtoutfile{--text follows this line--}
\immediate\write\gtoutfile{Proxy-for: \ifx\theasciiauthors\relax
\theauthors\else\theasciiauthors\fi\s<\ifx\theasciiemail\relax\theemail\else\theasciiemail\fi>}
\immediate\write\gtoutfile{\noexpand\\}
\immediate\write\gtoutfile{Authors: \ifx\theasciiauthors\relax
\theauthors\else\theasciiauthors\fi}
{\def\\{ }\immediate\write\gtoutfile{Title: \ifx\theasciititle\relax
\thetitle\else\theasciititle\fi}}
\immediate\write\gtoutfile{Subj-class: GT or SG, GR etc}
\immediate\write\gtoutfile{MSC-class: \theprimaryclass\ifx\thesecondaryclass\relax\else, \thesecondaryclass\fi}
\immediate\write\gtoutfile{Journal-ref: Algebraic and Geometric Topology \thevolumenumber\s
(\thevolumeyear) \startpage-\finishpage}
\immediate\write\gtoutfile{Comments: Published by Algebraic and
Geometric Topology at}
\immediate\write\gtoutfile{\s\s\s  http://www.maths.warwick.ac.uk/agt/AGTVol\thevolumenumber/agt-\thevolumenumber-\thepapernumber.abs.html}
\immediate\write\gtoutfile{\noexpand\\}
\immediate\write\gtoutfile{}
\ifx\theasciiabstract\relax
\immediate\write\gtoutfile{\theabstract}\else
\immediate\write\gtoutfile{\theasciiabstract}\fi
\immediate\write\gtoutfile{}
\immediate\write\gtoutfile{\noexpand\\}
\immediate\write\gtoutfile{}
\immediate\closeout\gtoutfile}}  

\def\maketitlepage{\makeagttitle\makeheadfile}
\let\makeshorttitle\maketitlepage
\let\maketitle\maketitlepage

\lognumber{3}
\volumenumber{1}
\volumeyear{2001}
\papernumber{3}
\published{9 December 2000}
\pagenumbers{39}{55}
\received{18 October 2000}\accepted{30 November 2000}

\usepackage{amssymb} 
\usepackage{amsmath,amscd} 
\usepackage{graphics} 

\newtheorem{theorem}{Theorem}[section]
\newtheorem{lemma}[theorem]{Lemma}
\newtheorem{cor}[theorem]{Corollary}
\newtheorem{proposition}[theorem]{Proposition}
\theoremstyle{definition}

\newtheorem{problem}[theorem]{Problem}
\theoremstyle{remark}
\newtheorem{remark}[theorem]{Remark}
\numberwithin{equation}{section}
\newcommand{\Q}{\mathbb Q}
\newcommand{\C}{\mathbb C}
\newcommand{\D}{\mathcal{D}}
\newcommand{\F}{\mathcal{F}}

\newcommand{\G}{\Gamma}
\newcommand{\gr}{\mathcal{G}}
\newcommand{\Gd}{\gr \delta}
\newcommand{\HE}{\mathcal H}
\newcommand{\K}{\mathcal K}
\newcommand{\h}{\mathfrak h}
\newcommand{\I}{{\mathcal I}}
\newcommand{\sgn}{\operatorname{sgn}}

\newcommand{\M}{{\mathcal M}}

\newcommand{\p}{\psi}
\newcommand{\x}{\xi}
\newcommand{\z}{\zeta}
\newcommand{\Z}{{\mathbb Z}}
\newcommand{\End}{\operatorname{End}}

\renewcommand{\ker}{\operatorname{Ker}}
\newcommand{\im}{\operatorname{Im}}
\newcommand{\var}[1]{\varphi^{(#1)}}
\begin{document}

\title[The Jones representation]{An expansion of the Jones 
representation\\of genus $2$ and the Torelli group}
\author{Yasushi Kasahara}
\address{Department of Electronic and Photonic System Engineering, 
  Kochi University of Technology, Tosayamada-cho, Kagami-gun, Kochi, 
  782--8502 Japan}
\email{kasahara@ele.kochi-tech.ac.jp}

\begin{abstract} 

We study the algebraic property of the representation of the mapping
class group of a closed oriented surface of genus $2$ constructed by
V\,F\,R Jones \cite{jones}. It arises from the Iwahori--Hecke algebra
representations of Artin's braid group of $6$ strings, and is defined
over integral Laurent polynomials ${\Z}[t, t^{-1}]$.  We substitute
the parameter $t$ with $-e^{h}$, and then expand the powers $e^h$ in
their Taylor series. This expansion naturally induces a filtration on
the Torelli group which is coarser than its lower central series. We
present some results on the structure of the associated graded
quotients, which include that the second Johnson homomorphism factors
through the representation. As an application, we also discuss the
relation with the Casson invariant of homology $3$--spheres.
\end{abstract}

\asciiabstract{We study the algebraic property of the representation
of the mapping class group of a closed oriented surface of genus 2
constructed by VFR Jones [Annals of Math. 126 (1987) 335-388]. It
arises from the Iwahori-Hecke algebra representations of Artin's braid
group of 6 strings, and is defined over integral Laurent polynomials
Z[t, t^{-1}].  We substitute the parameter t with -e^{h}, and
then expand the powers e^h in their Taylor series. This expansion
naturally induces a filtration on the Torelli group which is coarser
than its lower central series. We present some results on the
structure of the associated graded quotients, which include that the
second Johnson homomorphism factors through the representation. As an
application, we also discuss the relation with the Casson invariant of
homology 3-spheres.}



\makeshorttitle  

\section{Introduction}
\par

Let $\Sigma_{2}$ be a closed oriented surface of genus $2$, and let 
$\M_{2}$ be its mapping class group.
In \cite{jones}, V\,F\,R Jones constructed a finite dimensional linear 
representation of $\M_{2}$. We call this the Jones representation of 
genus $2$, and denote it by $\rho$  (see 
section \ref{def_of_rho} for the precise definition).
The Jones representation $\rho$ is defined over integral 
Laurent polynomials $\Z[t, t^{-1}]$. The construction of $\rho$ is based on the 
Iwahori--Hecke algebra representations of Artin's braid group of $6$ strings 
as well as the Birman--Hilden presentation of $\M_{2}$. The representation 
$\rho$ has been 
given attention as a candidate for the faithful linear representation 
of $\M_{2}$ (\cite{birman}). However, it seems that there are few 
results concerning $\rho$, with exceptions of 
a work of Humphries \cite{humphries} 
and the author's previous result \cite{kasahara}, 
both of which deal with specializations at roots of unity.
\par

Let $\I_{2}$ denote the Torelli group of genus $2$. It is defined as the kernel of 
the classical symplectic representation $\M_{2} \to Sp(4, \Z)$ which 
is induced by the natural action of $\M_{2}$ on the first integral 
homology group $H_{1} (\Sigma_{2}; \Z)$ of $\Sigma_{2}$. 
Here, $Sp(4, \Z)$ denotes the Siegel modular group.  In view of the 
explicit computation of $\rho$ given by Jones, it is easy to see that there are 
several special values of $t$ at which $\rho$ becomes trivial on $\I_{2}$. 
Such specializations include $t=\pm 1$.
It is natural to consider perturbations of $t$ at such special 
values to  extract any information of $\I_{2}$ in $\rho$.
\par

In this paper, we consider the perturbation at $t=-1$.
More precisely, let $\varphi$ be the expansion of $\rho$
obtained by setting $t = -e^{h}$ and then expanding the powers 
$e^h$ in their Taylor series. 
We then introduce a descending filtration of the Torelli group 
$\I_{2}$ by using $\varphi$,  and present some results on 
the associated graded quotients.  Our results should be considered as new 
properties of $\rho$. Our approach is based on  Jones' explicit computation
mentioned above.
\par

Now we state our main results. We set 
$$ \F_{k} = \I_{2} \cap \{ f \in \M_{2}; \, \varphi(f) \equiv 
\text{identity modulo terms of degrees higher than $k$} \}. $$
It will turn out that $\F_{1}$ coincides with $\I_{2}$ so that we have a  
filtration
$$ \I_{2}= \F_{1} \supset \F_{2} \supset \cdots \supset \F_{k} \supset \cdots .$$
A general result by Lazard in  \cite{lazard} can be used to deduce that
(a) for each $k \geq 1$, the  graded quotient 
$\gr_{k}\F = \F_{k}/\F_{k+1}$ is a free $\Z$--module of finite rank;
(b) the  filtration $\{ \F_{k} \}$ is central so that the associated graded 
sum $\gr\F = \oplus_{k \geq 1} \gr_{k}\F$  naturally forms a graded 
Lie algebra over $\Z$.  
Furthermore, it turns out that 
$\gr_{k}\F_{\Q} = \gr_{k}\F \otimes \Q$ has a natural structure of a 
rational $Sp(4, \Q)$--module where $Sp$ denotes the symplectic group.
This enables us to use classical symplectic representation theory to 
study the structure of $\gr\F$. 
\par

Let $\G_{a,b}$ denote the irreducible rational representation 
of $Sp(4, \Q)$ with highest weight $aL_{1} + b ( L_{1} + L_{2} )$ 
(cf section \ref{irrep}). We understand that the same symbol also 
denotes the representation space over $\Q$ in the obvious manner. 
Our first result is:
\begin{theorem} \label{A}
    Let $k=1$. Then the following isomorphism of $Sp(4, \Q)$--modules holds:
    $$ \gr_{1}\F_{\Q} = \I_{2} / \F_{2} \otimes \Q= \G_{0,2}$$
\end{theorem}
The appearance of the module $\G_{0,2}$ is suggestive since it also 
appears as the image of the second Johnson homomorphism $\tau_{2}(2)$ 
which is defined on $\I_{2}$ (see \cite{morita:survey} as well as for
the definition). In fact, we can show:
\begin{cor} \label{cor}
    The kernel of the second Johnson homomorphism $\tau_{2}(2)$ coincides 
    with $\F_{2}$. In particular, $\tau_{2}(2)$ factors through the Jones 
    representation $\rho$.
\end{cor}
Recall that Morita \cite{morita:casson} has described the Casson invariant 
of homology $3$--spheres in terms of mapping class groups of surfaces. 
In particular, he has related the second Johnson homomorphism with the 
Casson invariant. In view of this, Corollary \ref{cor} raises a problem whether 
the Jones representation $\rho$ contains the full information of the 
Casson invariant in the case of genus $2$. This 
problem will be reduced to the detection of a certain homomorphism of 
$\I_{2}$ which arises from Morita's theory of secondary 
characteristic classes of surface bundles 
\cite{morita:JDG}, \cite{morita:survey}.
We will discuss this later in section \ref{relation}.
\par

As for the module $\gr_{k}\F_{\Q}$ for general $k$, we have the following.
\begin{theorem} \label{B}
    For each $k \geq 2$, the following holds:
    \begin{enumerate}
          \item[\rm{(1)}] If $k$ is even, then $\gr_{k}\F_{\Q} \supset \G_{2,0}$.
          \item[\rm{(2)}] If $k$ is odd, then $\gr_{k}\F_{\Q} \supset \G_{0,2}$.
    \end{enumerate}
\end{theorem}
\par

Our expansion $t=-e^h$ was motivated by the fact that a similar 
expansion decomposes the Jones 
polynomial of links, which can be also defined by the Iwahori--Hecke 
algebra representations of braid groups,  into a series of Vassiliev 
invariants (see \cite{bar-natan}, \cite{birman-lin}, \cite{birman:bull}).  
The reader may consider that an expansion 
$t=e^h$ would be more natural. In this case, however, the 
graded quotients associated with the induced filtration have 
structures of $Sp(4, \Z/2\Z)$--modules, rather than $Sp(4, \Q)$--modules. 
Thus the induced filtration does not seem to fit the context of the Torelli group.
This phenomenon is due to the well known fact that the Iwahori--Hecke algebras 
(of type A) serve as deformations of group algebras of symmetric 
groups. We would like to study this aspect in detail elsewhere.

\par

The organization of the present paper is as follows. 
In section \ref{preliminary}, we recall some basic material concerned 
with the Jones representation $\rho$. We also determine the specialization 
of $\rho$ at $t=-1$, which will be the degree $0$ part of the expanded 
representation $\varphi$. 
In section \ref{filtration}, we describe the filtration $\{ \F_{k} \}$ in a 
more convenient manner and 
establish its basic properties.  In 
particular, we will see that each graded quotient $\gr_{k}\F$ is
naturally embedded into an $Sp(4, \Z)$--module which will be denoted by 
$\gr_{k}R$. In section \ref{str. of gr_k{R}}, we show that $\gr_{k}R$ is 
actually a rational $Sp(4, \Q)$--module, and then
we apply the symplectic representation theory to 
decompose $\gr_{k}R$ into weight spaces. In section \ref{proofs}, 
we complete the proofs of our main results.
Finally in section \ref{relation}, we discuss a possible relation with the 
Casson invariant.
\section{Preliminaries}
\label{preliminary}
\par

\subsection{Jones' construction} 
\begin{figure}[h]
   \begin{center}
   \includegraphics{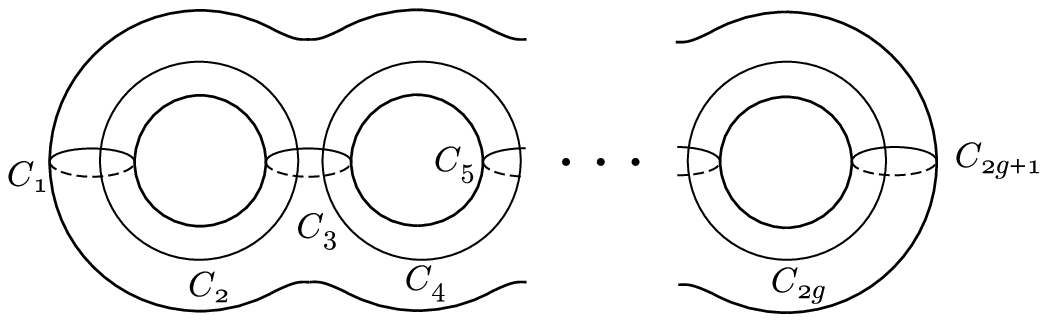}
      \nocolon
      \caption{}
      \label{fig.1}
   \end{center}
\end{figure}

We start with recalling the general 
construction by Jones concerning $\rho$. We refer to \cite{jones} 
for further details. Let $\Sigma_{g}$ be a closed oriented surface of genus $g$ and 
let $\M_{g}$ be its mapping class group. The hyperelliptic 
mapping class group $\HE_{g}$ is defined as the subgroup of 
$\M_{g}$ consisting of those elements which commute with the
class of a fixed hyperelliptic involution.  For $i=1$, $2$, \ldots, 
$(2g+1)$, let $\z_{i}$ be the Dehn twist along the simple closed curve $C_{i}$ 
on $\Sigma_{g}$ depicted in Figure \ref{fig.1}. 
Under an obvious choice of the 
hyperelliptic involution, each $\z_{i}$ lies in $\HE_{g}$. According 
to the presentation of $\HE_{g}$ due to Birman and Hilden 
\cite{birman-hilden}, $\HE_{g}$ is generated by $\z_{1}$, $\z_{2}$, \ldots, 
$\z_{2g+1}$, and these generators satisfy the defining relations of 
$B_{2g+2}$, Artin's braid group of $(2g+2)$ strings. We thus have a natural 
projection $p \co B_{2g+2} \to \HE_{g}$.  Recall that the irreducible 
representations of $B_{2g+2}$ which arise from the Iwahori--Hecke 
algebra of type $A_{2g+1}$ are in one-to-one correspondence with the 
Young diagrams of $(2g+2)$ boxes \cite{jones}. Jones considered which 
of such representations of $B_{2g+2}$ can yield a representation of 
$\HE_{g}$ via $p$. Let $Y$ be a Young diagram of $(2g+2)$ boxes, and 
let $\pi_{Y}$ denote the corresponding representation of 
$B_{2g+2}$. Jones proved that $\pi_{Y}$ defines, via $p$, a 
{\em projective} linear representation of $\HE_{g}$ if and only if $Y$ 
is {\em rectangular}. Furthermore, he proved that every projective 
representation of $\HE_{g}$ so obtained lifts to an actual linear 
representation of $\HE_{g}$ by determining the correcting scalar 
factor explicitly. Let us call this type of representation  the Jones 
representation of genus $g$. 
For example the $1 \times (2g+2)$ rectangular Young diagram gives the 
one dimensional trivial representation, and the $(2g+2) \times 1$ one gives 
a one dimensional representation given by the correspondence 
$\z_{i} \mapsto (-1) \cdot \text{identity}$. For later use, we denote the latter 
representation by $\sgn$. We will also denote suitable scalar extensions of 
$\sgn$ by the same symbol.
\par

\subsection{The column-row symmetry}
\par

For a Young diagram $Y$ with $(2g+2)$ boxes, let $Y'$ denote the 
Young diagram obtained by interchanging rows and columns in $Y$.  There 
exists an algebra involution on the Iwahori--Hecke algebra of type $A_{2g+1}$ 
which describes the transition from  each $Y$ to $Y'$ at the 
representation level (see \cite[Note 4.6]{jones}). We can combine this 
involution with Jones' construction above to obtain the relation 
between the two Jones representations corresponding to 
rectangular Young diagrams $Y$ and $Y'$, denoted by $\rho_{Y}$ and 
$\rho_{Y'}$, respectively: 
$$\rho_{Y'} (\z_{i}) = -\rho_{Y} (\z_{i}^{-1})$$
for $i=1$,\ldots, $(2g+1)$. Note that the correspondence 
$\z_{i} \mapsto \z_{i}^{-1}$ defines an automorphism $\iota$ of 
$\HE_{g}$, which can be realized by  the conjugation by an 
orientation {\em reversing} involution on $\Sigma_{g}$. Thus we have an equality 
of representations:
\begin{equation} \label{symmetry}
    \rho_{Y'} = (\rho_{Y} \circ \iota ) \otimes \sgn
\end{equation}
for every {\em rectangular} Young diagram $Y$ with $(2g+2)$ boxes.
\par

\subsection{The genus $2$ case} \label{def_of_rho}
\par

In general, $\HE_{g}$ is a proper subgroup of $\M_{g}$ and it is not 
obvious whether Jones representations of genus $g$ extend to $\M_{g}$ 
or not. However, in the case of $g=2$, it 
is classically known that $\HE_{2} = \M_{2}$ so that all the Jones 
representations of genus $2$ are defined on $\M_{2}$ for a trivial 
reason. Furthermore, in 
view of \eqref{symmetry},  the non-trivial Jones representation of 
genus $2$ is essentially unique.   This unique representation, which 
corresponds to  the $3 \times 2$ rectangular Young diagram, is the main object 
of this paper, and is denoted by $\rho$.  
\par

As was mentioned above, $\M_{2}$ is generated by $\z_{1}$, \ldots, $\z_{5}$. 
Jones has explicitly computed the images of these generators under 
$\rho$, which can be taken as the definition of $\rho$. Instead 
of $\z_{i}$'s, we use another set of generators of $\M_{2}$ 
which is more convenient for our 
computation. Let us take an element $\xi = \z_{1} \z_{2} \cdots 
\z_{5} \in \M_{2}$  which is a periodic automorphism of order $6$.
It is easy to see that $\z_{i+1} = \xi \z_{i} \xi^{-1}$ for 
$i = 1$, $2$,\ldots, $4$. Thus we can choose $\z_{1}$ and $\xi$ to be 
generators of $\M_{2}$. Now, due to the computation by Jones, 
the Jones representation $\rho \co  \M_{2} \to GL(5, \Z[t, t^{-1}])$ can 
be  defined 
by the following:
\begin{align*}
{\rho}(\z_1) = &
\begin{pmatrix}
 -1/t^2   &  0           &  0         &  0         &  t^3 \\
 0         &  -1/t^2     &  1/t^2    &  0         &  0   \\   
 0         &  0           & t^3        &  0         &  0   \\
 0         &  0           & 1/t^2     &  -1/t^2   &  0   \\
 0         &  0           &  0         &  0         &  t^3
\end{pmatrix}, \\
\intertext{and}
{\rho}(\x)= &
\begin{pmatrix}
 0 &  0 &  1 &  0 &  0 \\
 0 &  0 &  0 &  0 &  1 \\
 1 &  0 &  0 &  0 &  0 \\
 0 &  1 &  0 &  0 &  0 \\
 0 &  0 &  0 &  1 &  0
 \end{pmatrix}.
\end{align*}
Here our parameter $t$ corresponds to the formal power 
$q^{1/5}$ in \cite{jones}.
\par

\begin{figure}[h]
 \begin{center}
   \includegraphics{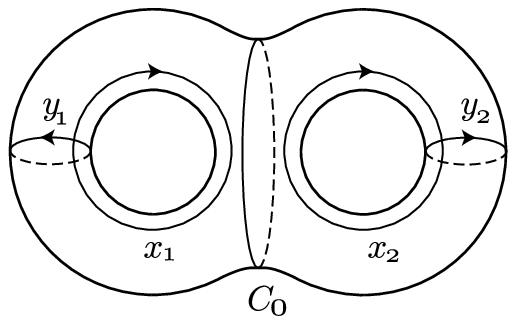}
   \nocolon
   \caption{}
   \label{fig.2}
 \end{center}
\end{figure}

Recall that the Torelli group of genus $2$, denoted by $\I_{2}$, 
is the kernel of the symplectic representation 
$\M_{2} \to Sp(4, \Z)$. Due to 
Powell \cite{powell}, $\I_{2}$ is the normal closure in $\M_{2}$ of 
the single Dehn twist $\p_{0}$ along the separating simple closed 
curve $C_{0}$ on $\Sigma_{2}$ depicted in Figure \ref{fig.2}. 
The image of $\p_{0}$ 
under $\rho$ has also been computed by Jones. In our notation, 
\begin{multline} \label{normal-generator}
    \rho(\p_{0})  = \rho(\z_{1} \z_{2} \z_{1})^4  =  t^6 \cdot Id  \\ + 
    \frac{t^{15}+1}{t^{24}} \cdot 
    \begin{pmatrix}
             0 & 0 & 0 & 0 & 0 \\
             0 & 0 & 0 & 0 & 0 \\
             0 & 0 & 0 & 0 & 0 \\ 
             t^{10}-1 &  -t^{10}+t^{5} & t^{10}-1 & -t^{15}+1 & -t^{10}+t^{5} \\
             0 & 0 & 0 & 0 & 0 
     \end{pmatrix}
\end{multline}
where $Id$ denotes the identity matrix. In particular, this shows that 
$\rho$ does not factor through the 
symplectic representation. On the other hand, it is easy to see that 
there exist several special values of $\rho$ at which $\rho$ becomes 
trivial on $\I_{2}$. Such special values include $t = \pm 1$. 
\par

We now determine the specialization of $\rho$ at $t=-1$. This will play a 
fundamental role in our investigation as the degree $0$ part of our expansion 
of $\rho$. Let $H$ denote the first integral homology group $H_{1} (\Sigma_{2}; \Z)$ of 
$\Sigma_{2}$, and let $\Lambda^2 H$ denote the second exterior product of $H$. 
The symplectic class $\omega \in \Lambda^2 H$ is defined as 
$x_{1} \wedge y_{1} + x_{2} \wedge y_{2}$ for any symplectic basis $x_{i}$, 
$y_{i}$ ($i=1$, $2$) of $H$ with respect to the algebraic intersection pairing.
 \begin{lemma} \label{pre-1}
     The specialization of $\rho$ at $t=-1$ descends to a linear 
     representation of $Sp(4,\Z)$ via the symplectic representation. 
     Furthermore, this specialization  is equivalent to 
     $(\Lambda^2 H / \omega \cdot \Z) \otimes \sgn$. 
\end{lemma}
 For later use, we fix a symplectic basis of $H$ 
and a basis of $\Lambda^2 H / \omega \cdot \Z$. Let $x_{1}$, $x_{2}$, 
$y_{1}$, $y_{2}$ be the symplectic basis of $H$ depicted in Figure 
\ref{fig.2}. We choose a basis $\{t_{i} ; i=1, \ldots, 5 \}$ of 
$\Lambda^2 H / \omega \cdot \Z$ as follows:
\begin{align} \label{basis}
    t_{1}& =[x_{1} \wedge x_{2}], & t_{2}& = [ y_{1} \wedge y_{2} ], &
    t_{3}& = [ x_{1} \wedge y_{1} ], \\
    t_{4}& = [ x_{1} \wedge y_{2} ], &
    t_{5}& = [ x_{2} \wedge y_{1} ]   \notag
\end{align}
where $[ \phantom{x}]$ denotes the corresponding class in
$\Lambda^2 H / \omega  \Z$.

\begin{proof}[Proof of Lemma \ref{pre-1}]
    As was mentioned above, the former part of the lemma
    is a direct consequence of \eqref{normal-generator}. 
    For the latter part, let $P$ denote the specialization of $\rho$ obtained 
    by putting $t=-1$.  Note that $P$ is defined over $\Z$ so that its 
    representation space is $\Z^5$. Let $Z$ and $X$ denote the matrix form of 
    the action of $\z_{1}$ and $\xi$, respectively,  on $\Lambda^2 H / \omega  \Z$ 
    with respect to the basis $\{ t_{i} \}$ above.  Then the assertion 
    follows from the existence of $F \in GL(5, \Z)$ such that 
    $P(\z_{1})F=-FZ$ and $P(\xi)F =-FX$. A direct calculation shows 
    that the following matrix gives a solution:
    \begin{equation} \label{F}
	F=
	\begin{pmatrix}
	    0 & -1 & 0 & 0 & 0 \\ 
	    0 & 0  & 0 & 0 & -1 \\ 
	    1 & 0  & 0 & 0 & 0 \\ 
	    0 & 0  & 1 & 1 & -1 \\ 
	    0 & 0  & 0 & 1 & 0 
	\end{pmatrix}
     \end{equation}
\vglue-.7cm
\end{proof}
\par
\section{The filtration $\{ \F_{k} \}$} \label{filtration}
\par
\subsection{The description of $\F_{k}$} 
\par
As in Introduction, we obtain the representation 
$$ \varphi \co \M_{2} \to R^{\times}$$
from $\rho$ by putting $t=-e^h$ and 
then expanding the powers $e^h$ in their Taylor series. 
Here $R$ denotes the algebra $M(5, \Q[[h]])$ of $5 \times 5$ matrices 
over rational formal power series in a variable $h$,  and $^\times$ 
denotes the group of invertible elements of an algebra. Hereafter, the 
identity of an algebra will be denoted by $1$. 
\par

For $k \geq 0$, let $R[k]$ be the principal ideal of $R$ generated by 
$h^{k} \cdot 1$.  The ideal $R[k]$ consists of those elements in $R$ 
the minimal degrees of whose terms are at least $k$. Thus we have a filtration 
of $R$:
$$ R=R[0] \supset R[1] \supset R[2] \supset \cdots $$
Note that $R[k] \cdot R[l] = R[k+l]$.  
Let $\gr_{k} R$ be the graded quotient $R[k] / R[k+1]$ and let 
$\gr R = \oplus_{k \geq 0} \gr_{k} R$ be the associated graded algebra.
The standard bracket $[X, Y] = XY - YX$ in $R$ induces a structure of 
a graded Lie algebra over $\Q$ on $\gr R$.
\par

Let  $\varphi^{(k)} \co \M_{2} \to (R/R[k+1])^{\times}$ be the composition of 
$\varphi$ with the natural projection $R^{\times} \to (R/R[k+1])^{\times}$. 
Now the filtration $\{ \F_{k} \}$ can be described as 
$\F_{k} = \I_{2} \cap \ker{\varphi^{(k-1)}}$ for $k \geq 1$.  Obviously, 
$\varphi^{(0)}$ is equivalent to the specialization of $\rho$ at $t=-1$. 
Thus we have $\F_{1} = \I_{2}$ by Lemma \ref{pre-1}.  As before, we 
write $\gr_{k} \F$ for $\F_{k} / \F_{k+1}$, and $\gr \F$ for 
$\oplus_{k \geq 1} \gr_{k} \F$. 
\par

\subsection{The graded Lie algebra $\gr \F$}
\par

We define a mapping $\delta \co \I_{2} \to R$ by 
$\delta(\p) = \varphi(\p) - 1$. It is obvious that $\delta(\F_{k}) \subset R[k]$.
For each $k \geq 1$, let $\delta_{k}$ be the composition of $\delta$
restricted to $\F_{k}$ with the obvious projection $R[k] \to \gr_{k}R$. 
The following proposition can be deduced from a general result of 
Lazard \cite{lazard} (see also \cite{MKS}).  
\begin{proposition} \label{main-1}
    Each $\delta_{k}$ induces an embedding of the module  $\gr_{k}\F$ 
    into $\gr_{k}R$. Furthermore, the commutator operation 
    $[ x, y ] = x y x^{-1} y^{-1}$ induces a  structure of a graded Lie algebra on 
    $\gr \F$, and $\oplus \delta_k$ gives an embedding 
    of a graded Lie algebra $\gr \F \to \gr R$.
\end{proposition}
\begin{remark} \label{main-2}
    The quotient module $\gr_{k}R$ has a natural structure of a finite 
    dimensional vector space over $\Q$. Thus the proposition  
    implies that each $\gr_{k}\F$ is a free abelian group of finite rank. 
\end{remark}
We give a direct proof of Proposition \ref{main-1} for completeness.
\begin{lemma} \label{main-3}
   For $k \geq 1$, the mapping 
   $\delta_{k} \co \F_{k} \to \gr_{k}R =  R[k] / R[k+1]$ is a group 
   homomorphism.
\end{lemma}
\proof
    If $\p_{1}$, $\p_{2} \in \F_{k}$, then we have $\delta(\p_{1}) 
    \delta(\p_{2}) \in R[2k] \subset R[k+1]$. Thus the lemma follows 
    immediately from the following equality:
    $$ \delta(\p_{1} \p_{2}) = \delta(\p_{1}) + 
                                      \delta(\p_{2}) + \delta(\p_{1}) \delta(\p_{2}). \eqno{\qed}$$

 Clearly $\ker{\delta_{k}} = \F_{k+1}$. Thus $\delta_{k}$ 
gives an embedding of  a group 
$$ \gr_{k} \F \to \gr_{k} R,$$ 
which will be denoted by the same symbol $\delta_{k}$.
\par
Now let $\gr \delta$ denote $\oplus_{k \geq 1} \delta_{k} \co \gr\F 
\to \gr R$. To prove the latter part of the proposition, we have 
to see that $\im{\gr \delta}$ is a Lie subalgebra of $\gr R$, and 
that the Lie bracket on $\gr \F$ induced by $\gr \delta$ corresponds to the
commutator in $\I_{2}$.  These follow from the next lemma.
\begin{lemma} \label{main-4}
   Let $x \in \F_{k}$ and $y \in \F_{l}$. Then $ [ x, y ] \in 
   \F_{k+l}$, and the following equality holds:
   $$ \delta_{k+l} ( [ x, y ] ) = \delta_{k} (x)  \delta_{l} (y)  -  
   \delta_{l} (y) \delta_{k} (x)$$
\end{lemma}
\begin{proof}
    A direct computation implies the following equality in $R$:
    \begin{equation} \label{main-eq:1}
         \delta( [ x, y ] ) = (\delta(x) \delta(y) - \delta(y) 
         \delta(x) ) \varphi( x^{-1} y^{-1} )
    \end{equation}
    This implies immediately that $ [ x, y ] \in \F_{k+l}$. 
    Furthermore, in view of the fact that  $\varphi^{(0)}$, the degree $0$ part of 
    $\varphi$, is trivial on $\I_{2}$ (Lemma \ref{pre-1}), we obtain the 
    required equality by taking \eqref{main-eq:1} modulo $\F_{k+l+1}$.
\end{proof}
This completes the proof of Proposition \ref{main-1}.
\subsection{$Sp(4, \Z)$--actions} \label{symp_action}
\par
For each $k$, $\gr_{k}\F$  naturally has a structure of an 
$Sp(4, \Z)$--module in the following manner. The Torelli group 
$\I_{2}=\F_{1}$ is acted on by $\M_{2}$ via conjugation. The 
filtration $\{ \F_{k} \}$ is preserved by this action of $\M_{2}$ so 
that $\gr_{k}\F$ is  an $\M_{2}$--module. Since 
$[ \I_{2}, \F_{k} ] \subset \F_{k+1}$ by Lemma \ref{main-4}, the 
action of $\I_{2}$ on $\gr_{k}\F$ is trivial. Hence the above 
action of $\M_{2}$ descends to that of $Sp(4, \Z)$.
\par

Next we consider the natural action of $Sp(4, \Z)$ on $\gr_{k}R$ for 
each $k \geq 0$. Recall that $\var{k}$
denotes the composition of the expanded Jones representation $\varphi$ with 
the projection $R^{\times} \to (R/R[k+1])^{\times}$. For each $\zeta \in \M_{2}$, 
the correspondence 
$x \mapsto \zeta_{*}x = \var{k}(\zeta)  x  \var{k}(\zeta^{-1})$ induces a 
structure of an $\M_{2}$--module on $R/R[k+1]$. Thus $\gr_{k}R \subset 
R/R[k+1]$ is an $\M_{2}$--submodule. It is easy to see that the following 
formula holds
\begin{equation}
    \zeta_{*}x = \var{0}(\zeta) \cdot x \cdot \var{0}(\zeta^{-1})
    \label{pf-eq:2}
\end{equation}
for $\zeta \in \M_{2}$ and $x \in \gr_{k}R$. Here the multiplication 
in the right hand side is meant by the one in  the  graded algebra $\gr R$.
Now the triviality of $\var{0}$ on $\I_{2}$ (Lemma \ref{pre-1}) implies 
that the action of $\M_{2}$ on $\gr_{k}R$ factors through 
$Sp(4, \Z)$.
\par

\begin{proposition} \label{pf-1}
    The embedding of Lie algebra  $\Gd \co \gr\F \to \gr R$ is 
    $Sp(4, \Z)$--equivariant. 
\end{proposition}

\begin{proof}
    It suffices to show that 
    $ \delta_{k}(\zeta \p \zeta^{-1}) = \zeta_{*} \delta_{k}(\p)$
    for $\zeta \in \M_{2}$, $\p \in \F_{k}$. By definition, we have 
    $\delta( \zeta \p \zeta^{-1} ) = \varphi(\zeta) \delta(\p) 
    \varphi(\zeta^{-1})$. Thus taking the mod $R[k+1]$ classes, we have
    $\delta_{k}(\zeta \p \zeta^{-1}) = \var{0}(\zeta) \cdot \delta_{k}(\p) \cdot 
    \var{0}(\zeta^{-1})$ as  elements of $\gr R$. In view of 
    \eqref{pf-eq:2}, this is nothing but $\zeta_{*} \delta_{k}(\p)$.
\end{proof}
\par
\section{The structure of $\gr_{k}R$} \label{str. of gr_k{R}}
\par

\subsection{Irreducible representations of  $Sp(4, \Q)$} \label{irrep}
\par
Now we recall the description of irreducible 
representations of the algebraic group $Sp(4, \Q)$. We follow the book 
of Fulton and Harris \cite{FH}. Let ${\mathfrak{sp}}(4, \C)$ be the Lie 
algebra of the symplectic Lie group $Sp(4, \C)$, and let ${\h}$ 
be its Cartan subalgebra consisting of diagonal matrices. Choose a 
system of fundamental weights $L_{1}$ and $L_{2} \co {\h} \to \C$ 
as in \cite{FH}. Then for each pair $(a, b)$ of non-negative integers, 
there exists a unique irreducible representation of $Sp(4, \C)$ with 
highest weight $a L_{1} + b (L_{1} + L_{2})$. We denote this 
representation by $\Gamma_{a, b}$ following \cite{FH}. These 
are all rational representations defined over $\Q$ so that we can 
consider them as irreducible representations of $Sp(4, \Q)$.  We will 
understand that the same notation $\Gamma_{a, b}$ also denotes the 
representation space over $\Q$ of the corresponding representation in 
the obvious manner. For example the trivial representation $\Q = 
\Gamma_{0, 0}$, $H_{\Q} = H_{1} (\Sigma_{2}; \Q) = \Gamma_{1, 0}$
and $\Lambda^2 H_{\Q} / \omega  \Q = \Gamma_{0, 1}$ where $\omega$ 
denotes the symplectic class. 
\par
\subsection{The action of $Sp(4,\Q)$ on $\gr_{k}R$}
\par

For an $Sp(4, \Z)$--module $V$, let $V_{\Q}$ denote $V \otimes \Q$.
In view of  \eqref{pf-eq:2}, it is easy to see that the $Sp(4, \Z)$--module 
$\gr_{k}R$ is isomorphic to $(\var{0} \otimes (\var{0})^{*})_{\Q}$ where 
$*$ denotes the dual of a representation. 
By Lemma \ref{pre-1}, we have
$\var{0}_{\Q} = \G_{0,1} \otimes \sgn$,  which is only an 
$Sp(4, \Z)$--module.
However, the triviality of $\sgn \otimes \sgn^{*} = \sgn \otimes \sgn$ implies 
that 
$$
    (\var{0} \otimes (\var{0})^{*})_{\Q} = \G_{0,1} \otimes 
    \G^{*}_{0,1} = \G_{0,1} \otimes \G_{0,1}
$$
so that $\gr_{k}R$ can be naturally considered as an $Sp(4,\Q)$--module. 
Here we have used the self-duality for general symplectic modules to deduce the last 
equality. Furthermore, the decomposition  of  $\G_{0,1} \otimes \G_{0,1}$
into irreducible modules is well known (eg  \cite[section 16.2]{FH}), 
which implies the equality
\begin{equation} \label{pf-eq:3}
    \gr_{k}R = \G_{0,2} \oplus \G_{2,0} \oplus \G_{0,0}.
\end{equation}
\subsection{Weight spaces of $\gr_{k}R$}
\par

Now we explicitly describe the decomposition of
$\gr_{k}R = \G_{0,1} \otimes \G_{0,1}^{*}$ into weight spaces and 
determine the submodule $\G_{0,2}$  in terms of weight spaces. 
Recall that we have fixed the basis $\{ t_{i}; i=1, \ldots, 5 \}$ of 
$\Lambda^2 H / \omega \Z$ in \eqref{basis}. The $\{t_{i} \}$ also 
serves as a basis of $\G_{0,1} = (\Lambda^2 H / \omega \Z) \otimes \Q$ 
in the obvious manner. All the members $t_{i}$'s are weight 
vectors in $\G_{0,1}$ and the corresponding weights are: $(L_{1} + L_{2})$, 
$-(L_{1} + L_{2})$,  $0$, $(L_{1} - L_{2})$ and $(-L_{1} + L_{2})$ 
for $i=1$, \ldots, $5$, respectively. These can be deduced from the 
description of the weight spaces for the fundamental representation 
$H_{\Q}=H_{1}(\Sigma_{2}; \Q)$.
\par
Now let $\{ t_{i}^{*} \}$ denote the dual basis of $\{ t_{i} \}$, and 
let $e_{i,j}$ denote $t_{i} \otimes t_{j}^{*} \in \G_{0,1} \otimes 
\G_{0,1}^{*}$. Then $\{ e_{i,j} \}$ forms a basis of $\G_{0,1} 
\otimes \G_{0,1}^{*}$. Furthermore,  each $e_{i,j}$ is a weight 
vector and its weight is given by  the weight of $t_{i}$ minus that of 
$t_{j}$. Consequently, we can obtain the decomposition of $\G_{0,1} 
\otimes \G_{0,1}^{*}$ into weight spaces as given by Table 
\ref{pf:eq-3.5}.  Each row of Table \ref{pf:eq-3.5} consists of a 
weight appearing in $\G_{0,1} \otimes \G_{0,1}^{*}$ and a basis for
the corresponding weight space.
\begin{table}[ht!] 
    \caption{The weight spaces of $\G_{0,1} \otimes \G_{0,1}^{*}$}
    \label{pf:eq-3.5}
    \begin{center}\small
	\begin{tabular}{r|l}
	    \noalign{\hrule height0.8pt}
	    \hfil Weight & A basis for the weight space  \\
	    \hline
	    $2(L_{1} + L_{2})$ & $\{e_{1,2} \}$ \\
	    $L_{1} + L_{2}$ & $\{e_{1,3}, e_{3,2} \}$ \\
	    $2 L_{1}$ & $\{ e_{1,5}, e_{4,2} \}$ \\
	    $2L_{2}$ & $\{e_{1,4}, e_{5,2} \}$ \\
	    $2(L_{1} - L_{2})$ & $\{ e_{4,5} \}$ \\
	    $L_{1} - L_{2}$ & $\{e_{3,5}, e_{4,3} \}$ \\
	    $0$ & $\{ e_{i,i};  i=1, \ldots, 5 \}$ \\
	    $-L_{1} + L_{2}$ & $\{ e_{3,4}, e_{5,3} \}$ \\
	    $2(-L_{1} + L_{2} )$ & $\{ e_{5,4} \}$ \\
	    $-2L_{2}$ & $\{ e_{2,5}, e_{4,1} \}$ \\
	    $-2L_{1}$ & $\{ e_{2,4}, e_{5,1} \}$ \\
	    $-(L_{1} + L_{2})$ & $\{ e_{2,3}, e_{3,1} \}$ \\
	    $-2(L_{1} + L_{2})$ & $\{e_{2,1} \}$ \\
	 \noalign{\hrule height0.8pt}
	 \end{tabular}
     \end{center}
\end{table}
The identification of  $\G_{0,1} \otimes \G_{0,1}^{*} = 
\End{(\G_{0,1})}$ with $\gr_{k}R$ is given as follows. Let 
$E_{i,j} \in M(5, \Q)$ denote the matrix form of $e_{i,j} \in 
\End{(\G_{0,1})}$ with respect to the basis $\{ t_{i} \}$. 
As was observed in the proof of Lemma \ref{pre-1}, $\G_{0,1}$ 
is identified with $\var{0} \otimes \sgn$ via the matrix $F$ given by 
\eqref{F}. 
Then the desired identification is given by the correspondence
\begin{equation} \label{pf-eq:4}
    e_{i,j} \mapsto h^{k} F E_{i,j} F^{-1} \mod{R[k+1]}.
\end{equation}
\par

Next we describe the weight spaces of $\G_{0,2}$. Note that 
$\G_{0,2}$ is the highest weight submodule of 
$\G_{0,1} \otimes \G_{0,1}^{*}$. 
Recall from the general representation theory that the highest weight 
submodule is generated by the images of a highest weight vector,
which are also weight vectors, under 
successive applications of  the negative root spaces of the Lie algebra. 
Due to this fact, we can obtain the weight spaces  of $\G_{0,2}$ starting 
with a highest weight vector $e_{1,2}$. The 
result is given by Table \ref{pf:eq-5}.
\begin{table}[ht!] 
    \caption{The weight spaces of $\G_{0,2}$}
    \label{pf:eq-5}
    \begin{center}\small
	\begin{tabular}{r|l}
	    \noalign{\hrule height0.8pt}
	    \hfil Weight & A basis for the weight space  \\
	    \hline
	    $2(L_{1} + L_{2})$ & $\{ e_{1,2} \}$ \\
	    $L_{1} + L_{2}$ & $\{ e_{1,3} + 2e_{3,2} \}$ \\
	    $2L_{1}$ & $\{ e_{1,5} + e_{4,2} \}$ \\
	    $2L_{2}$ & $\{ e_{1,4} + e_{5,2} \}$ \\
	    $2(L_{1} - L_{2} )$ & $\{ e_{4,5} \}$ \\
	    $L_{1} - L_{2}$ & $\{ e_{4,3} + 2e_{3,5} \}$ \\
	    $0$ & $ \{ (e_{1,1} + e_{2,2}) - (e_{4,4} + e_{5,5}), 
	                     2e_{3,3} - ( e_{1,1} + e_{2,2}) \}$ \\
	    $L_{2} - L_{1}$ & $\{ 2e_{3,4} + e_{5,3} \}$ \\
	    $2(L_{2} - L_{1})$ & $\{ e_{5,4} \}$ \\
	    $-2L_{2}$ & $\{ e_{2,5} + e_{4,1} \}$ \\
	    $-2L_{1}$ & $\{ e_{2,4} + e_{5,1} \}$ \\
	    $-(L_{1} + L_{2})$ & $\{ e_{2,3} + 2e_{3,1} \}$ \\
	    $-2(L_{1} + L_{2})$ & $\{ e_{2,1} \}$ \\
	    \noalign{\hrule height0.8pt}
	 \end{tabular}
     \end{center}
\end{table}
\par
\section{Proofs of Theorems} \label{proofs}
\par

First of all, by virtue of the following general proposition, the action of 
$Sp(4, \Z)$ on $\gr_{k}\F$ described in section \ref{symp_action} 
can be extended to an action of 
$Sp(4,\Q)$ on $\gr_{k}\F_{\Q} = \gr_{k}\F \otimes \Q$ so that $\Gd \otimes \Q$ is an 
$Sp(4, \Q)$--equivariant homomorphism 
between graded Lie algebras over $\Q$.
\begin{proposition}[Asada--Nakamura {\cite[(2.2.8)]{asada-nakamura}}]  
    \label{pf-2}
    Let $g \geq 1$.
    If $L$ is a $\Z$--submodule of a rational finite dimensional 
    $Sp(2g, \Q)$--module $V$ which is stable under the action of 
    $Sp(2g, \Z)$, then $L \otimes \Q \subset V$ is also stable under 
    the action of $Sp(2g, \Q)$.
\end{proposition}
\par

\subsection{Proof of Theorem \ref{A}}
\par
In view of \eqref{pf-eq:3},
it suffices to show that $\delta_{1}(\F_{1})_{\Q} = \G_{0,2} 
\subset \gr_{1}R$. As was explained in section \ref{def_of_rho}, 
$\F_{1} = \I_{2}$ is the normal closure of $\p_{0}$, the Dehn twist along 
the separating simple closed curve $C_{0}$. Thus $\delta_{1} (\F_{1})$ is the minimal 
$\Z$--submodule of $\gr_{1}R$ which contains $\delta_{1}(\p_{0})$ 
and is invariant under the action of $Sp(4,\Z)$. A direct computation shows 
that 
$$ \delta_{1}(\p_{0}) \equiv 
     F \cdot \begin{pmatrix}
	6h & 0  &     0  &   0 & 0 \\ 
	 0 & 6h &     0  &   0 & 0 \\ 
	 0 & 0   & -24h &   0 & 0 \\ 
	 0 & 0   &     0  & 6h & 0 \\ 
	 0 & 0   &     0  &   0 & 6h 
      \end{pmatrix} \cdot F^{-1} 
    \mod{R[2]}
$$
where the matrix $F$ is the one given in  \eqref{F}.
This corresponds, under the identification \eqref{pf-eq:4}, to 
$$-12 ( 2e_{3,3} - ( e_{1,1} + e_{2,2} ) ) 
      - 6 ( ( e_{1,1} + e_{2,2} ) - ( e_{4,4} + e_{5,5} ) )  \in 
      \G_{0,1} \otimes \G_{0,1}^{*}.$$
In view of Table \ref{pf:eq-5}, this is a weight vector of $\G_{0,2} 
\subset \G_{0,1} \otimes \G_{0,1}^{*}$ with weight $0$. Thus we have
$\delta_{1}(\F_{1}) \subset \G_{0,2}$. Now we can apply Proposition 
\ref{pf-2} to conclude that 
$\delta_{1}(\F_{1})_{\Q} =  \G_{0,2}$.
\qed
\par

\subsection{Proof of Corollary \ref{cor}}
\par

We freely use the notation of \cite[section 5]{morita:survey} for Johnson 
homomorphisms. To prove Corollary \ref{cor}, we have only to prove 
that there exists an $Sp(4, \Q)$--equivariant isomorphism 
$f \co \delta_{1}(\I_{2})_{\Q} \to \tau_{2}(2)(\I_{2})_{\Q} = 
\tau_{2}(2)(\I_{2}) \otimes \Q$ such that 
\begin{equation} \label{pf-eq:6}
    f \circ \delta_{1}(\p_{0})=\tau_{2}(2)(\p_{0}).
\end{equation}
As mentioned in Introduction, it is well known that $\tau_{2}(2)(\I_{2})_{\Q}$ 
is isomorphic to $\G_{0,2}$, and hence there exists an essentially unique 
$Sp(4, \Q)$--isomorphism  between $\delta_{1}(\I_{2})_{\Q}$ 
and $\tau_{2}(2)(\I_{2})_{\Q}$. However, the existence of $f$ 
satisfying the condition \eqref{pf-eq:6} does not follow from this 
fact directly. We need an analysis of the images of $\p_{0}$ under 
both $\delta_{1}$ and $\tau_{2}(2)$ in terms of highest weight vectors.
\par

Recall that $\delta_{1}(\I_{2})_{\Q}$ is 
generated by the images of its highest weight vector under the action 
of the universal enveloping algebra $U$ of $\mathfrak{sp}(4,\Q)$. 
Thus $\delta_{1}(\p_{0})$ can be expressed as $Z \cdot e_{1,2}$ for 
some fixed element $Z \in U$.
\par

On the other hand, an explicit computation of $\tau_{g,1}(2)$, a 
relative version of $\tau_{g}(2)$, has been 
given by Morita \cite{morita:casson}. We can combine this result with
the relationship between $\tau_{g,1}(2)$ and $\tau_{g}(2)$ explained 
in \cite{morita:survey} to obtain the description of $\tau_{g}(2)$.  
More precisely, we can describe (a) the decomposition of 
$\tau_{2}(2)(\I_{2})_{\Q}$ into weight spaces together with the action 
of $\mathfrak{sp}(4,\Q)$; (b) an expression of $\tau_{2}(2)(\p_{0})$ 
in terms of weight vectors.  As a result, we can choose a highest 
weight vector $v_{0}$ of $\tau_{2}(2)(\I_{2})_{\Q}$ such that 
$\tau_{2}(2)(\p_{0}) = Z \cdot v_{0}$. If we take $f$ to be the unique 
$Sp(4, \Q)$--equivariant isomorphism which sends $e_{1,2}$ to $v_{0}$, 
then it satisfies the condition \eqref{pf-eq:6}. This completes the 
proof of Corollary \ref{cor}.
\qed
\par

\subsection{Proof of Theorem \ref{B}}
\par

We show that the desired submodules come from the lower central series 
of the Torelli group $\{ \I [k] \}$ where $\I[1] = \I_{2}$ and 
$\I[k+1]$ is defined recursively  by  $\I[k+1]  =  [ \I_{2}, \I[k] \, ]$ 
for $k \geq 1$.  More precisely, let 
$\gr \D_{\Q}$ denote the Lie subalgebra of $\gr \F_{\Q}$ generated 
by $\gr_{1}\F_{\Q}$. We write $\gr_{k}\D_{\Q}$ for its homogeneous 
component of degree $k$ so that we have $\gr_{1}\D_{\Q} = 
\gr_{1}\F_{\Q}$ and $\gr_{k+1}\D_{\Q} = [\gr_{1}\D_{\Q}, \gr_{k}\D_{\Q} ]
\subset \gr_{k+1}\F_{\Q}$ for $k \geq 1$. It is obvious that 
$\gr_{k}\D_{\Q} = \delta_{k}(\I [k] / \I [k+1]) \otimes \Q$.  We assert that 
$$\gr_{k}\D_{\Q} = 
\begin{cases}
    \G_{2,0}  \text{ for } k \text{ even} \\
    \G_{0,2}  \text{ for } k \text{ odd}
\end{cases}
$$
as a submodule of $\End{(\G_{0,1})}$ under the identification of 
\eqref{pf-eq:4}. Under the same identification, the Lie 
bracket $[\; ,\: ] \co \gr_{k}R \otimes \gr_{l}R \to \gr_{k+l}R$ corresponds 
to the standard Lie bracket in $\End{(\G_{0,1})}$ given by 
$[A,B] = AB-BA$. Since we have already seen that 
$\gr_{1}\D_{\Q} = \G_{0,2}$, the assertion above follows from the 
next two equalities in $\End{(\G_{0,1})}$:
\begin{align*}
    & [\G_{0,2}, \G_{0,2}] = \G_{2,0}\\
    & [\G_{0,2}, \G_{2,0}] = \G_{0,2}
\end{align*}
These equalities can be checked as follows. First we can easily see 
that the Lie bracket is $Sp(4,\Q)$--equivariant so that both the 
left hand sides above are direct sums of some of $\G_{0,2}$, 
$\G_{2,0}$ and $\G_{0,0}$.
\par

Next, we can use Table \ref{pf:eq-5} to compute the dimensions of weight 
spaces  in $[\G_{0,2}, \G_{0,2}]$ with weights $2(L_{1} + L_{2})$, 
$2L_{1}$ and $0$, respectively.  The results are $0$, $1$ and $2$, 
respectively. Thus the highest weight of $[\G_{0,2}, \G_{0,2}]$ 
is $2L_{1}$, which implies that 
$\G_{0,2} \not\subset [\G_{0,2}, \G_{0,2}]$  and  
$\G_{2,0} \subset [\G_{0,2}, \G_{0,2}]$. Finally we can see that 
$\G_{0,0} \not\subset [\G_{0,2}, \G_{0,2}]$ by comparing the 
dimension of the weight space of weight $0$ in $[\G_{0,2}, \G_{0,2}]$
with that in $\G_{2,0}$ (cf \cite{FH}). This concludes the first equality.
\par

Now with the assistance of the first equality, we can easily describe the 
whole weight spaces of $\G_{2,0} = [\G_{0,2}, \G_{0,2}]$. Then the 
second equality follows from a similar argument.
This completes the proof of Theorem \ref{B}.
\qed
\section{The relation with the Casson invariant} \label{relation}
\par
In this section, motivated by Corollary \ref{cor}, we discuss a 
possible relation between the Jones representation $\rho$ and the 
Casson invariant of homology $3$--spheres, and pose a related 
problem.

\par
We first recall the description of the Casson invariant, denoted by 
$\lambda$, in terms of mapping class groups of surfaces by Morita. We 
refer to \cite{morita:casson}, \cite{morita:torelli} for further 
details. As before, let $\Sigma_{g}$ be a closed oriented surface of 
genus $g \geq 2$, and let $\M_{g}$ be its mapping class group. Also, 
let $\K_{g}$ be the subgroup of $\M_{g}$ generated by all the Dehn 
twists along {\em separating} simple closed curves on $\Sigma_{g}$. Given a 
homology $3$--sphere $M$ together with an embedding $f \co \Sigma_{g} \to M$,  
the Casson invariant gives rise to a homomorphism 
$\lambda_{f} \co \K_{g} \to \Z$. Roughly, $\lambda_{f}$ is defined by 
$\lambda_{f}(\psi) = \lambda( M_{\psi} ) - \lambda ( M )$ for 
$\psi \in \K_{g}$ where $M_{\psi}$ is the homology $3$--sphere
obtained from $M$ by first cutting along $f ( \Sigma_{g} )$ and then 
regluing by $\psi$. Furthermore, $\lambda_{f}$ can be 
decomposed as a sum of two $\Q$--valued homomorphisms:
$$ \lambda_{f} = -\frac{1}{24} d_{1}  +  q_{f}. $$
Here the homomorphism $d_{1}$ coincides with a secondary invariant 
which arises from Morita's theory of secondary characteristic classes of surface 
bundles (\cite{morita:casson}, \cite{morita:JDG}, \cite{morita:survey}).  In particular, 
$d_{1}$ is independent of the choice of the pair $(M, f)$, and is 
invariant under the conjugation action of $\M_{g}$ on $\K_{g}$.  It is 
known that such a $\Q$--valued homomorphism of $\K_{g}$ is unique up to 
nonzero scalars (\cite[Theorem 5.7]{morita:torelli}). 
\par

On the other hand, the homomorphism $q_{f}$ factors through 
$\tau_{g}(2)$ and depends on the choice of $(M, f)$.

\par

Now we assume that $g=2$. Note that $\K_{2} = \I_{2}$ as 
mentioned before. In view of Corollary \ref{cor},  we see that 
$q_{f}$ factors through $\rho$. Thus we may say that $\rho$ contains 
the full information of $\lambda$, in the case of $g=2$, if and only 
if $\rho$ contains the information of $d_{1}$.
\par

Now we set $\mathcal{Z} = \I_{2} / \ker{d_{1}}$, and 
$\mathcal{Z}_{\sigma} = \sigma ( \I_{2} ) / \sigma ( \ker{d_{1}} )$ 
for an arbitrary homomorphism $\sigma$ of $\I_{2}$. Note that 
$\mathcal{Z}$ is isomorphic to the  infinite cyclic group $\Z$ and 
$\mathcal{Z}_{\sigma}$ is a quotient of $\mathcal{Z}$. 
Under these notation,  the information of $d_{1}$ contained in $\rho$ 
concentrates on the cyclic group $\mathcal{Z_{\rho}}$. On the other 
hand, we can easily see that $\ker{d_{1}} = [ \M_{2}, \I_{2} ]$.  We 
thus have the following:
\begin{proposition} \label{rel-prop}
    The homomorphism $d_{1}$ factors through the Jones 
    representation $\rho$ if and only if 
    \begin{equation} \label{rel-1}
	\mathcal{Z}_{\rho} = \rho( \I_{2} ) / \rho( [ \M_{2}, \I_{2} ] )
	  \cong \Z.
    \end{equation}  
\end{proposition}
At the moment, we do not know whether the equality \eqref{rel-1} is 
true or not.
\begin{problem}
    Determine the order of the cyclic group $\mathcal{Z}_{\rho}$.
\end{problem}
\begin{remark}
  \begin{enumerate}
      \item[\rm{(1)}]  If the order of $\mathcal{Z}_{\rho}$ is finite, then it 
          implies that $\rho$ is {\em not faithful} on $\I_{2}$.
       \item[\rm{(2)}] Our previous result \cite{kasahara} implies that the order of 
          $\mathcal{Z}_{\rho}$ is at least $10$.
       \item[\rm{(3)}] For any reduction $r$ of $\rho$, the computation of 
           $\mathcal{Z}_{r}$ would be helpful in estimating the order of 
           $\mathcal{Z}_{\rho}$.  With such an expectation, we considered 
           the cases $r = \varphi^{(1)}$ and $r=\varphi^{(2)}$. As a result, we 
           obtained the isomorphisms 
	   $$ \mathcal{Z}_{\varphi^{(1)}} = \mathcal{Z}_{\tau_{2}(2)} \cong 
	    \Z/10,  \qquad \mathcal{Z}_{\varphi^{(2)}}  \cong \Z/10
	    $$
            by computer calculations.
    \end{enumerate}	    
\end{remark}
\par

In view of the above, it remains possible that the order of 
$\mathcal{Z}_{\rho}$ is exactly $10$. In any case, it might be interesting 
to note that $\Z/10$ is isomorphic to the abelianization of $\M_{2}$.

%
%

\Addresses
\recd

\end{document}